\documentclass[11pt]{article}

\usepackage[a4paper,margin=32mm]{geometry}
\usepackage[T1]{fontenc}
\usepackage{lmodern}
\usepackage{amsmath,amssymb,amsthm,mathtools}
\usepackage{booktabs}
\usepackage{microtype}
\usepackage{enumitem}
\usepackage{xcolor}
\definecolor{linkblue}{RGB}{25,65,125}
\usepackage{tabularx}
\usepackage{hyperref}
\usepackage[nameinlink,capitalise,noabbrev]{cleveref}

\hypersetup{
  pdftitle={A Geometric Characteristic Long Exact Sequence for Codimension-Two Spin and Pin-plus/minus Bordism},
  pdfauthor={Kohei Takeda},
  colorlinks=true,
  linkcolor=linkblue,
  citecolor=linkblue,
  urlcolor=linkblue
}

\newtheorem{theorem}{Theorem}[section]
\newtheorem*{introtheorem}{Theorem}
\newtheorem{lemma}[theorem]{Lemma}
\newtheorem{proposition}[theorem]{Proposition}
\newtheorem{corollary}[theorem]{Corollary}
\theoremstyle{definition}
\newtheorem{definition}[theorem]{Definition}
\newtheorem{remark}[theorem]{Remark}

\newcommand{\Ftwo}{\mathbb Z/2}

\newcommand{\Spin}{\mathrm{Spin}}
\newcommand{\Pinp}{\mathrm{Pin}^{+}}
\newcommand{\Pinm}{\mathrm{Pin}^{-}}
\newcommand{\Pinpm}{\mathrm{Pin}^{\pm}}
\newcommand{\PD}{\operatorname{PD}}
\newcommand{\Th}{\operatorname{Th}}

\newcommand{\id}{\operatorname{id}}
\newcommand{\intt}{\operatorname{int}}

\title{A Geometric Characteristic Long Exact Sequence\\
for Codimension-Two Spin and $\Pinpm$ Bordism}
\author{Kohei Takeda}
\date{}

\begin{document}

\maketitle

\begin{abstract}
Following Botvinnik--Rosenberg \cite{BotvinnikRosenberg2023}, we
define fibered-boundary bordism, introduce sphere-bundle bordism,
and prove an exact sequence relating them to ordinary bordism.
The sequence also holds when all sphere fibers are required
to be non-$D^r$-bounding. Following Debray--Ye--Yu
\cite{DebrayYeYu2026}, we define characteristic cycles with a
chosen complement structure trivializing a characteristic class.
For the GM and $\mathrm{KT}^{\pm}$ cases, deleting tubular
neighborhoods and identifying sphere-bundle structures with
determinant-twisted base structures give a characteristic exact
sequence.  This proves the bordism-group
consequence of the DYY conjecture for our definitions, without
constructing the conjectured map of spectra.
\end{abstract}
\section{Introduction}

Bordism with fibered singularities studies manifolds whose
boundaries carry a specified fiber-bundle structure.
Botvinnik--Rosenberg \cite{BotvinnikRosenberg2023} develop this
approach in their work on positive scalar curvature.
Following their framework, we define a fibered-boundary bordism
theory for sphere bundles of real vector bundles and arbitrary
stable tangential structures.  The purpose of this paper is to
relate this theory to characteristic bordism, where an embedded
submanifold represents a characteristic class of the ambient
manifold.  We construct geometric
isomorphisms between the relevant bordism groups and use them to
derive a characteristic long exact sequence.

\paragraph{Bordism with fibered boundary.}
Let $\psi:B\psi\to BO$ be a stable tangential structure and fix
$r\geq1$.  A degree-$n$ fibered-boundary cycle is a tuple
$(X,B,p,V,\varphi)$, where $X$ is a compact $n$-dimensional
$\psi$-manifold, $B$ is a closed $(n-r)$-manifold,
$p:\partial X\to B$ is a smooth $S^{r-1}$-bundle, and $V\to B$
is a real rank-$r$ vector bundle.  The map
$\varphi:S(V)\to\partial X$ is a sphere-bundle isomorphism over
$B$, so $p\circ\varphi=\pi_V$.  Pullback by $\varphi$ gives $S(V)$ the boundary structure of $X$.  Bordisms extend these data on
their fibered faces, as specified in \cref{def:fib}.
We write the resulting group as $\Omega_n^{\psi\text{-}\mathrm{fib}}$.

We also introduce the group $\Omega_k^{\psi\text{-}\mathrm{sph}}$ of
sphere-bundle cycles $(B,V,\eta)$, where $B$ is a closed
$k$-manifold, $V$ has rank $r$, and $\eta$ is a $\psi$-structure
on $TS(V)$.  A bordism consists of a bundle over a compact base
bordism, with a structure on its sphere-bundle total space
extending the structures at the ends; see \cref{def:sph}.
No separate $\psi$-structure on the base is required.
There are three maps: $i$ views a closed manifold as a cycle
with empty fibered boundary, $\beta$ takes the fibered boundary
with its bundle data, and $\tau$ takes the total space $S(V)$.

\begin{introtheorem}[\ref{thm:fib}]
For any stable tangential structure $\psi$ and any $r\geq1$,
these maps form a long exact sequence
\[
 \cdots\longrightarrow\Omega_{n-r+1}^{\psi\text{-}\mathrm{sph}}
 \xrightarrow{\tau}\Omega_n^\psi
 \xrightarrow{i}\Omega_n^{\psi\text{-}\mathrm{fib}}
 \xrightarrow{\beta}\Omega_{n-r}^{\psi\text{-}\mathrm{sph}}
 \xrightarrow{\tau}\Omega_{n-1}^\psi\longrightarrow\cdots.
\]
\end{introtheorem}
The proof follows the gluing argument of
\cite[Theorem~3.15]{BotvinnikRosenberg2023}, using the sphere-bundle
bordism group to record the boundary data.

\paragraph{The nonbounding theories.}
For a rank-$r$ bundle $V\to B$, a $\psi$-structure on $TS(V)$
is $D^r$-bounding on a fiber if its restriction, with the
constant base directions retained, is induced from a structure
on $TD(V_b)\oplus\underline{T_bB}$ by the boundary convention.
We define $\Omega_*^{\mathrm{nb},\psi\text{-}\mathrm{fib}}$ and
$\Omega_*^{\mathrm{nb},\psi\text{-}\mathrm{sph}}$ by requiring
non-$D^r$-bounding fibers on both cycles and bordisms.
The exact sequence remains valid for any stable $\psi$ and
any $r\geq1$ (\cref{prop:nb}).

\paragraph{Characteristic cycles.}
Characteristic bordism studies manifolds with submanifolds
representing Poincar\'e duals of characteristic classes.
Codimension-two examples appear in the work of Freedman--Kirby
\cite{FreedmanKirby1978}, Guillou--Marin \cite{GuillouMarin1986},
and Kirby--Taylor \cite{KirbyTaylor1990}.
Debray--Ye--Yu \cite[Definition~2.26]{DebrayYeYu2026} give a
general definition of characteristic pairs $(M,F)$ for a
tangential structure $\xi$ and a characteristic class $P$.

For a stable tangential structure $\xi:B\xi\to BO$ and a class
represented by $P:B\xi\to K(\Ftwo,r)$, we define a
$\psi$-structure as a $\xi$-structure together with a chosen
trivialization of $P$.  Our characteristic cycle is a triple
$(M,F,\eta)$, where $M$ is a closed $n$-dimensional $\xi$-manifold,
$F\subset M$ is a closed embedded codimension-$r$ submanifold,
and $\eta$ is a chosen $\psi$-structure on $M\setminus F$ refining
the given $\xi$-structure and not extending across any component
of $F$.  Bordisms carry the same data, with a neat defect and a
complement structure that does not extend across any of its
components (\cref{def:char}).  We denote the group by
$\Omega_n^{(\xi,P)}$.

After this definition, we restrict to $r=2$ and
$(\xi,P,\psi)=(SO,w_2,\Spin)$, $(O,w_2+w_1^2,\Pinm)$, and
$(O,w_2,\Pinp)$ for the rest of the characteristic theory.
These are the GM, $\mathrm{KT}^-$, and $\mathrm{KT}^+$ cases,
respectively.  We interpret the chosen complement structure as
implicit in DYY's later twisted Spin descriptions and make it
explicit in our definition; see \cref{sec:related-work}.

The extension criterion of Kirby--Taylor
\cite[Theorem~2.4]{KirbyTaylor1990}, proved here in
\cref{thm:extension}, states that such a complement structure
exists if and only if
\[
 [F]=\PD_M(P(TM)).
\]
This explains the link with characteristic pairs.  The criterion
is an existence statement: the cycle includes a choice of a
structure with the required nonextension property.

\paragraph{The comparison isomorphisms.}
For the three triples above, nonextension across a component
of $F$ is equivalent to the nonbounding condition on its normal
$S^1$-fibers.  Deleting a tubular neighborhood therefore defines
\[
 [M,F,\eta]\longmapsto
 [M\setminus\intt D(\nu_{F/M}),F,
   \pi_{\nu_{F/M}},\nu_{F/M},\id],
\]
where the complement carries the restriction of $\eta$ and its
boundary is identified with $S(\nu_{F/M})$.
\Cref{thm:charfib} proves that this map is an isomorphism
\[
 \Omega_n^{(\xi,P)}\cong
 \Omega_n^{\mathrm{nb},\psi\text{-}\mathrm{fib}}.
\]
Its inverse attaches $D(V)$ along $\varphi$, takes the zero
section as the characteristic submanifold, and extends the
complement structure fiberwise over the punctured disk bundle.

To describe the other term in the sequence, we use
$TS(V)\cong\pi^*(TB\oplus\det V)$.
\Cref{lem:circle} proves that pullback gives a bijection between
$\psi$-structures on $TB\oplus\det V$ and non-$D^2$-bounding
$\psi$-structures on $TS(V)$, compatible with bordisms.
It follows that
\[
 \Omega_k^{\mathrm{nb},\psi\text{-}\mathrm{sph}}
 \cong\Omega_k^\psi\bigl((BO_2)^{\sigma-1}\bigr),
\]
where $\sigma$ is the universal determinant line and
$(BO_2)^{\sigma-1}$ is its virtual rank-zero Thom spectrum;
our tangential convention is given in \cref{rem:thom}.
Combining these two identifications gives the main result.

\begin{introtheorem}[\ref{thm:main}]
For each of the three triples $(\xi,P,\psi)$ above, there is a
long exact sequence
\[
 \cdots\longrightarrow
 \Omega_{n-1}^{\psi}\bigl((BO_2)^{\sigma-1}\bigr)
 \xrightarrow{t}\Omega_n^\psi
 \xrightarrow{i}\Omega_n^{(\xi,P)}
 \xrightarrow{R}\Omega_{n-2}^{\psi}\bigl((BO_2)^{\sigma-1}\bigr)
 \xrightarrow{t}\Omega_{n-1}^\psi\longrightarrow\cdots.
\]
\end{introtheorem}
Here $R$ takes the characteristic submanifold with its normal
bundle and the structure on $TF\oplus\det\nu_{F/M}$ determined
by $\eta$.  The map $t$ takes the total space of the sphere bundle
with the corresponding $\psi$-structure.

Debray--Ye--Yu \cite[Conjecture~2.32]{DebrayYeYu2026} conjecture
that $R$ is induced by a map of spectra.  This would place $R$
in a long exact sequence of bordism groups
\cite[Equation~(2.33)]{DebrayYeYu2026}.
For the characteristic theories defined here, \cref{thm:main}
proves this bordism-group consequence without assuming the
conjecture, and identifies the remaining term as ordinary Spin,
$\Pinm$, or $\Pinp$ bordism.  We do not construct the conjectured
map of spectra or a spectrum-level fiber sequence.

The paper is organized as follows.  \Cref{sec:related-work}
reviews earlier work.  \Cref{sec:fibered} defines fibered-boundary
and sphere-bundle bordism, proves their exact sequence, and then
treats their nonbounding versions and exact sequence.
\Cref{sec:characteristic} defines characteristic bordism, restricts
to the three obstruction problems, and proves the two comparison
isomorphisms and the main theorem.

\section{Related work}\label{sec:related-work}

Freedman--Kirby \cite{FreedmanKirby1978} use oriented
characteristic surfaces in oriented four-manifolds in a geometric
proof of Rochlin's theorem.  Guillou--Marin
\cite{GuillouMarin1986} allow nonorientable characteristic surfaces
in an oriented four-manifold and prove a signature congruence.

Kirby--Taylor \cite{KirbyTaylor1990} describe the local geometry of Spin and $\Pinpm$
structures.  Their \S1 describes the $\Pinpm$ obstruction classes and
the $H^1(-;\Ftwo)$-torsors of lifting structures; Lemma~1.7 gives
determinant-line relations.  Theorem~2.4 is the codimension-two
extension criterion.  Their Corollary~6.12 gives a
Freedman--Kirby bordism exact sequence, and Remarks~6.15--6.16
discuss related constructions.  The ordinary Spin and $\Pinpm$ bordism
theories are also studied in
\cite{AndersonBrownPeterson1967,AndersonBrownPeterson1969,
KirbyTaylorPinPlus1990}.

These sources use different characteristic data.  Kirby--Taylor's \S6 includes data for
$w_1$ and compatibility conditions for an induced structure on the
defect.  Debray--Ye--Yu explain how their conventions differ
from these in \cite[Definition~2.31 and Remark~3.87]{DebrayYeYu2026}.

Debray--Ye--Yu \cite{DebrayYeYu2026} describe characteristic and
defect structures by twisted Spin theories.  For GM their defect
term is determinant-twisted Spin bordism
\cite[\S3.3]{DebrayYeYu2026}.  We compare their $\Pinpm$-type defect formulas
with ours by computing the obstruction classes.
Let $\sigma\to BO(1)$ and $V\to BO_2$ be the tautological
real line bundle and real $2$-plane bundle, respectively,
and use the same notation for their pullbacks to the submanifold $F$ in a
$\mathrm{KT}^\pm$ characteristic pair.
For $\mathrm{KT}^-$,
\cite[Lemma~3.64]{DebrayYeYu2026} gives
\[
w_1(TF)=w_1(\sigma)+w_1(V),\qquad
w_2(TF)=w_1(\sigma)^2+w_1(\sigma)w_1(V)+w_1(V)^2.
\]
Substituting $w_1(\sigma)=w_1(TF)+w_1(V)$ into the second
identity and applying the Whitney formulas gives
\[
w_2(TF\oplus\det V)+w_1(TF\oplus\det V)^2=0.
\]
For $\mathrm{KT}^+$,
\cite[Lemma~3.90]{DebrayYeYu2026} gives
$w_1(TF)=w_1(\sigma)+w_1(V)$ and
$w_2(TF)=w_1(V)^2+w_1(\sigma)w_1(V)=w_1(TF)w_1(V)$.
Hence $w_2(TF\oplus\det V)=0$.
We interpret the complement $\psi$-structure, not explicitly listed
in their characteristic cycles, as part of their later
twisted Spin descriptions.  With this interpretation, our
convention agrees with theirs.

Our use of fibered boundaries is related to bordism with prescribed
singularities \cite{Baas1973} and manifolds with fibered singularities
\cite{BotvinnikRosenberg2023}.  We prove exactness directly
for the given sphere bundles with the required fiber condition.
It does not use the low-dimensional characteristic bordism
calculations in these sources.

\section*{Acknowledgements}
I am deeply grateful to my supervisor, Takuya Sakasai, for his patient guidance, thoughtful advice, and constant encouragement throughout this project.

\paragraph{Use of AI.}
I used AI tools to proofread the manuscript, search for relevant
literature, check references, and discuss details of the proofs.
I take responsibility for the mathematical arguments and references
in this paper.

\section{Bordism with fibered boundary}\label{sec:fibered}

\subsection{Stable tangential structures and conventions}
\label{sec:structures}

Throughout this section, let $\psi:B\psi\to BO$ be a stable
tangential structure and fix an integer $r\geq1$.
A $\psi$-structure on a bundle $E\to X$ consists of a lift of its
stable classifying map to $B\psi$, together with the identifying
homotopy over $BO$.

We write $\psi(E)$ for the set of isomorphism classes of stable
$\psi$-structures on a fixed bundle $E$, where isomorphisms cover
the identity of $E$.  Diffeomorphisms of structured manifolds
include identifications of their structures.
We use product structures near boundaries.  The boundary
structure is defined by the outward-normal identification
\[
 TW|_{\partial W}\cong\mathbb R_{\mathrm{out}}\oplus T\partial W.
\]
For a $\psi$-manifold $X$, give $X\times[0,1]$ the product
structure whose induced structure on $X\times\{1\}$ is the
given structure on $X$.  We write $\overline X$ for $X$ with
the structure induced on $X\times\{0\}$, using the outward
normal $-\partial_t$.  Thus the two ends are
$\overline X\sqcup X$ with their induced structures.
For oriented bordism, $\overline X$ has the opposite orientation.
Nonnegative-dimensional cycles generate the groups, and groups of
negative degree are zero.

\subsection{Bordism with fibered boundary}

For a real rank-$r$ vector bundle $V\to B$, write $S(V)$ and $D(V)$
for its unit sphere and disk bundles.  Giving
$V$ up to bundle isomorphism is equivalent to giving its classifying
map $f:B\to BO(r)$ up to homotopy.  We use the bundle notation to
make the boundary identifications explicit.

\begin{definition}\label{def:fib}
A degree-$n$ $\psi$-fibered-boundary cycle is a tuple
$(X,B,p,V,\varphi)$ with the following data:
\begin{enumerate}[label=(\roman*),leftmargin=*]
\item $X$ is a compact smooth $n$-dimensional $\psi$-manifold;
\item $B$ is a closed $(n-r)$-manifold;
\item $p:\partial X\to B$ is a specified smooth $S^{r-1}$-fiber
bundle;
\item $V\to B$ is a real rank-$r$ vector bundle;
\item $\varphi:S(V)\xrightarrow{\cong}\partial X$ is an
isomorphism of smooth $S^{r-1}$-fiber bundles over $B$, so that
\begin{equation}\label{eq:bundle-phi}
 p\circ\varphi=\pi_V,
 \qquad \pi_V:S(V)\longrightarrow B.
\end{equation}
\end{enumerate}
The boundary $\psi$-structure induced by $X$ determines a
$\psi$-structure on $S(V)$ via $\varphi$.

Following the fibered-boundary bordism convention of
\cite[Definition~3.14]{BotvinnikRosenberg2023}, a bordism from
$(X_0,B_0,p_0,V_0,\varphi_0)$ to $(X_1,B_1,p_1,V_1,\varphi_1)$
consists of a compact $(n+1)$-dimensional $\psi$-manifold $W$
with corners and two boundary faces
\[
 \partial W=\partial_{\mathrm o}W\cup\partial_{\mathrm f}W,
 \qquad \partial_{\mathrm o}W\cong X_1\sqcup\overline X_0,
\]
a compact $(n-r+1)$-manifold $N$ with $\partial N=B_0\sqcup B_1$,
a smooth $S^{r-1}$-fiber bundle $q:\partial_{\mathrm f}W\to N$,
a real rank-$r$ bundle $\mathcal V\to N$, and a sphere-bundle
isomorphism
\begin{equation}\label{eq:bundle-Phi}
 \Phi:S(\mathcal V)\xrightarrow{\cong}\partial_{\mathrm f}W,
 \qquad q\circ\Phi=\pi_{\mathcal V}.
\end{equation}
At the common corner
$\partial(\partial_{\mathrm o}W)=\partial(\partial_{\mathrm f}W)
=q^{-1}(\partial N)$, the data $(q,\mathcal V,\Phi)$ restrict to
$(p_i,V_i,\varphi_i)$ over $B_i$.  All data are products near
the common corner.
We denote the bordism group by $\Omega_n^{\psi\text{-}\mathrm{fib}}$.
\end{definition}

\begin{definition}\label{def:sph}
A degree-$n$ $\psi$-sphere-bundle cycle is $(B,V,\eta)$, where
$B$ is a closed $n$-manifold, $V\to B$ is a real rank-$r$ vector bundle,
and $\eta$ is a $\psi$-structure on $TS(V)$, whose dimension
is $n+r-1$.
A bordism is $(N,\mathcal V,\eta_N)$, where $N$ is a compact
$(n+1)$-manifold, $\mathcal V$ restricts to the specified boundary
bundles, and the $\psi$-structure $\eta_N$ on $TS(\mathcal V)$
induces the given outgoing structure and the incoming structure
defined by the cylinder convention on
\[
 \partial S(\mathcal V)=S(\mathcal V|_{\partial N}).
\]
We do not require a separate $\psi$-structure on $TN$.
We denote the bordism group by $\Omega_n^{\psi\text{-}\mathrm{sph}}$.
\end{definition}

Isomorphisms of sphere-bundle cycles consist of a diffeomorphism
$g:B\to B'$ and a vector-bundle isomorphism $A:V\to V'$ covering
$g$, whose induced sphere map $S(A)$ preserves the specified
structure.  An isomorphism of fibered-boundary cycles also includes
a structured diffeomorphism $h:X\to X'$ and satisfies
\begin{equation}\label{eq:cycle-isomorphism}
 p'\circ\partial h=g\circ p,\qquad
 \partial h\circ\varphi=\varphi'\circ S(A).
\end{equation}
We use the chosen unit sphere bundles to define the sphere maps.
Bordisms include the same type of identifications at their ends,
preserving $q$, $\Phi$, and their boundary restrictions.

We keep $r$ fixed and omit it from the group notation.
We have maps
\[
\begin{aligned}
 i:\Omega_n^\psi&\longrightarrow\Omega_n^{\psi\text{-}\mathrm{fib}},
 & [M]&\longmapsto[M,\varnothing],\\
 \beta:\Omega_n^{\psi\text{-}\mathrm{fib}}&\longrightarrow
 \Omega_{n-r}^{\psi\text{-}\mathrm{sph}},
 &[X,B,p,V,\varphi]&\longmapsto[B,V,\varphi^*\eta_{\partial X}],\\
 \tau:\Omega_n^{\psi\text{-}\mathrm{sph}}&\longrightarrow\Omega_{n+r-1}^\psi,
 &[B,V,\eta]&\longmapsto[S(V),\eta].
\end{aligned}
\]
An ordinary bordism has empty fibered face, a fibered bordism
induces the sphere-bundle bordism
$(N,\mathcal V,\Phi^*\eta_{\partial_{\mathrm f}W})$ on its
fibered face, and a
sphere-bundle bordism has total space an ordinary bordism.
At the corner, exchanging the two outward normals reverses
the induced structure, so the fibered face gives a bordism
from $\beta(X_1)$ to $\beta(X_0)$, with the specified bundle data.
This still identifies their bordism classes.
Equation~\eqref{eq:cycle-isomorphism} shows that isomorphic cycles
have the same images.  Thus the three maps are well-defined.

\begin{theorem}[Fibered-boundary exact sequence]\label{thm:fib}
For any stable tangential structure $\psi$ and any $r\geq1$,
there is a long exact sequence
\[
 \cdots\longrightarrow\Omega_{n-r+1}^{\psi\text{-}\mathrm{sph}}
 \xrightarrow{\tau}\Omega_n^\psi
 \xrightarrow{i}\Omega_n^{\psi\text{-}\mathrm{fib}}
 \xrightarrow{\beta}\Omega_{n-r}^{\psi\text{-}\mathrm{sph}}
 \xrightarrow{\tau}\Omega_{n-1}^\psi\longrightarrow\cdots.
\]
\end{theorem}
\begin{proof}
We adapt the geometric proof of
\cite[Theorem~3.15]{BotvinnikRosenberg2023} to the sphere bundles
and tangential structures specified here.
First $\beta\circ i=0$ because the fibered boundary is empty, and
$\tau\circ\beta=0$ because the sphere bundle with its transported
structure is $\partial X$.  For $i\circ\tau=0$, use
$S(V)\times[0,1]$, with one end ordinary and the other end fibered
over $B$, with the specified projection $\pi_V$ and the identity
sphere-bundle identification at the fibered end.  This is a fibered null-bordism of $i([S(V),\eta])$,
with the incoming structure on its fibered end.

At the sphere-bundle term, let $[B,V,\eta]\in\ker\tau$.
Choose an ordinary null-bordism $X$ with structured boundary
$\partial X=S(V)$.  Taking $p=\pi_V$ and $\varphi=\mathrm{id}$
gives a fibered-boundary cycle whose image under $\beta$ is
$[B,V,\eta]$.

At the fibered-boundary term, suppose $\beta([X,B,p,V,\varphi])=0$.
Choose a sphere-bundle null-bordism $(N,\mathcal V,\eta_N)$ and
put $Y=S(\mathcal V)$.  Its boundary is identified with $S(V)$ as a sphere bundle over
$B$, with the prescribed structure.  Composing with $\varphi$
identifies it with $\partial X$ as a structured sphere bundle
with projection $p$.  Using that identification, glue
\[
 Z=X\cup_{\partial X}\overline Y.
\]
The boundary identifications allow the structures to glue, so
$Z$ is a closed $\psi$-manifold.  The gluing bordism is obtained
from $Z\times[0,1]$ by taking the upper end and the $X$ part of
the lower end as ordinary faces, and the $Y$ part as the fibered
face, introducing a corner along $\partial X$.
The latter face carries the projection and bundle identification
from $Y=S(\mathcal V)\to N$.  Their restrictions at the corner
are $p$ and $\varphi$, since $p\circ\varphi=\pi_V$.
This gives a bordism from $(X,B,p,V,\varphi)$ to $i(Z)$.

At the ordinary term, let $[M]\in\ker i$ and choose a fibered
null-bordism $W$ of $i(M)$.  Since $M$ is closed, its fibered face
has closed base $N$ and comes with $q$ and $\Phi$ from
\cref{def:fib}.  Use $\Phi$ to pull back its structure and obtain the
sphere-bundle cycle $(N,\mathcal V,\eta_{\mathrm f})$.
The ordinary boundary relation for $W$ gives
$[M]=\pm[S(\mathcal V),\eta_{\mathrm f}]$.  Thus $[M]$ lies
in $\operatorname{im}\tau$.  Together with the three zero composites,
these constructions prove exactness at every term.
\end{proof}

\subsection{The nonbounding fiber condition}
\label{sec:cycles}

\paragraph{The non-$D^r$-bounding condition on sphere fibers.}

Let $V\to B$ be a real vector bundle of rank $r$,
and let $\eta$ be a stable $\psi$-structure on $TS(V)$.
For each $b\in B$, the tangent bundle sequence gives a splitting $TS(V)|_{S(V_b)}
\cong TS(V_b)\oplus\underline{T_bB}$,
unique up to homotopy. Here $\underline{T_bB}$ denotes the
constant bundle with fiber $T_bB$.

We call $\eta$ \emph{$D^r$-bounding on the fiber over $b$}
if this restricted structure is induced, by the outward-normal
boundary convention of \cref{sec:structures}, from a stable
$\psi$-structure on $TD(V_b)\oplus\underline{T_bB}$. Otherwise, it is \emph{non-$D^r$-bounding on that fiber}.
We call $\eta$ \emph{fiberwise non-$D^r$-bounding}
if it is non-$D^r$-bounding on every fiber.

Homotopies of splittings and local trivializations, together
with products of disk extensions, show that this condition is
well-defined, locally constant on $B$, and unchanged by the
cylinder convention.

\begin{definition}\label{def:nb}
The group $\Omega_n^{\mathrm{nb},\psi\text{-}\mathrm{fib}}$
consists of the cycles of \cref{def:fib} for which the induced
structure on $S(V)$ is fiberwise non-$D^r$-bounding.
A bordism between two such cycles is a bordism
$(W,N,q,\mathcal V,\Phi)$ as in \cref{def:fib} for which
the induced structure on $S(\mathcal V)$ is also fiberwise
non-$D^r$-bounding.

Similarly, $\Omega_n^{\mathrm{nb},\psi\text{-}\mathrm{sph}}$
is defined using the sphere-bundle cycles and bordisms of
\cref{def:sph}, with the fiberwise non-$D^r$-bounding condition
imposed on both cycles and bordisms.
As before, the rank $r$ is fixed and omitted from the notation.
\end{definition}

\begin{proposition}\label{prop:nb}
For any stable tangential structure $\psi$ and any $r\geq 1$,
the maps $i$, $\beta$, and $\tau$ give a long exact sequence
\[
\cdots\longrightarrow
\Omega_{n-r+1}^{\mathrm{nb},\psi\text{-}\mathrm{sph}}
\xrightarrow{\tau}\Omega_n^\psi
\xrightarrow{i}\Omega_n^{\mathrm{nb},\psi\text{-}\mathrm{fib}}
\xrightarrow{\beta}
\Omega_{n-r}^{\mathrm{nb},\psi\text{-}\mathrm{sph}}
\xrightarrow{\tau}\Omega_{n-1}^\psi
\longrightarrow\cdots.
\]
\end{proposition}

\begin{proof}
The map $i$ uses an empty fibered boundary, and $\beta$ keeps
the given fiber structures. A sphere-bundle bordism has total
space an ordinary $\psi$-bordism, so $\tau$ is also well-defined.

The composites $\beta\circ i$ and $\tau\circ\beta$ vanish
for the same reasons as in \cref{thm:fib}.
The cylinder proving $i\circ\tau=0$ has a fibered end with
the incoming structure. The cylinder convention preserves
the non-$D^r$-bounding condition, so this is a bordism
in the theory of \cref{def:nb}.

At the sphere-bundle term, an element of $\ker\tau$ has an
ordinary filling whose boundary is the specified sphere-bundle
cycle. This filling is a fibered-boundary cycle with the
required fiber condition.

At the fibered-boundary term, an element of $\ker\beta$ has
a sphere-bundle null-bordism satisfying the fiber condition
throughout. The gluing construction in \cref{thm:fib} therefore
gives a bordism to a closed $\psi$-manifold whose fibered face
still satisfies that condition.

At the ordinary term, a fibered null-bordism of $i(M)$ has
a fibered face satisfying the condition by definition.
Its base is closed, so this face defines a sphere-bundle cycle.
The boundary relation places $[M]$ in the image of $\tau$.
These constructions prove exactness.
\end{proof}

\section{Characteristic bordism and the main exact sequence}
\label{sec:characteristic}

\subsection{Characteristic bordism as fibered-boundary bordism}
\label{sec:charcycles}

Let $r\geq1$, let $\xi:B\xi\to BO$ be a stable tangential
structure, and fix a map $P:B\xi\to K(\Ftwo,r)$ representing a
characteristic class.
A $\psi$-structure on a bundle $E\to X$ consists of a lift
$\ell:X\to B\xi$ of its stable classifying map, together with
the identifying homotopy over $BO$ and a chosen null-homotopy
$h:P\circ\ell\simeq *$.

\begin{definition}\label{def:char}
A degree-$n$ $(\xi,P)$-characteristic cycle is a triple $(M,F,\eta)$,
where $M$ is a closed smooth $n$-dimensional $\xi$-manifold,
$F\subset M$ is a closed embedded codimension-$r$ submanifold,
and $\eta$ is a chosen $\psi$-structure on $M\setminus F$ refining
the specified $\xi$-structure and not extending across any
component of $F$.

A bordism from $(M_0,F_0,\eta_0)$ to $(M_1,F_1,\eta_1)$
is a triple $(W,G,\eta_W)$, where $W$ is a compact
$(n+1)$-dimensional $\xi$-manifold, $G\subset W$ is a neat
codimension-$r$ submanifold, and $\eta_W$ is a $\psi$-structure
on $W\setminus G$ refining the given $\xi$-structure and
nonextendable across every component of $G$.
The boundary identifications are
\[
 (\partial W,\partial G)
 \cong (M_1,F_1)\sqcup(\overline M_0,F_0).
\]
They preserve the $\xi$-structures and identify the induced
complement structures with $\eta_1$ at the outgoing end and
with the cylinder-induced structure from $\eta_0$ at the
incoming end. Here $\partial G=G\cap\partial W$, and the
embedding and structures are products near $\partial W$.
Disjoint union defines the bordism group $\Omega_n^{(\xi,P)}$.
\end{definition}

From this point onward, set $r=2$ and fix one of the following
three triples:
\begin{equation}\label{eq:triples}
\begin{array}{c|c|c}
\xi&P&\psi\\ \hline
SO&w_2&\Spin\\
O&w_2+w_1^2&\Pinm\\
O&w_2&\Pinp
\end{array}
\end{equation}
For these three structures, the stabilized $S^1$-fiber has two
lifting types (with orientation fixed in the Spin case).
Exactly one extends over the disk; the other is
non-$D^2$-bounding in the sense of \cref{sec:cycles}.
For a normal bundle $\nu_{F/M}$, we use the boundary structure
of the complement of a tubular neighborhood.

\begin{theorem}[Characteristic extension criterion
{\cite[Theorem~2.4]{KirbyTaylor1990}}]\label{thm:extension}
Fix a triple $(\xi,P,\psi)$ in \eqref{eq:triples}.
Let $M$ be a closed smooth $n$-dimensional $\xi$-manifold and
$F\subset M$ a closed embedded codimension-two submanifold,
possibly disconnected.  There exists a $\psi$-structure on
$M\setminus F$ refining the specified $\xi$-structure and
nonextendable across every component of $F$ if and only if
\begin{equation}\label{eq:characteristic}
 [F]=\PD_M(P(TM))\quad\text{in }H_{n-2}(M;\Ftwo).
\end{equation}
\end{theorem}
\begin{proof}
We recall the argument of
\cite[Theorem~2.4]{KirbyTaylor1990}; see also
\cite[Proposition~2.24]{DebrayYeYu2026}.
Put $A=M\setminus F$. All cohomology groups below have
coefficients in $\Ftwo$.
The Thom isomorphism gives $H^2(M,A)\cong H^0(F)$.
Let $U_F$ correspond to the function equal to $1$ on every
component of $F$, and let $j:H^2(M,A)\to H^2(M)$ be the
natural map. Then $j(U_F)$ is Poincar\'e dual to $[F]$.

A chosen $\psi$-structure $\eta$ on $A$ has a relative
extension obstruction $o(\eta)\in H^2(M,A)$.
Indeed, the homotopy fiber of $B\psi\to B\xi$ is
$K(\Ftwo,1)$, so the obstruction lies in degree two and
there are no higher obstructions.
Under the Thom isomorphism, its value on each component
of $F$ is $0$ if $\eta$ extends over a normal disk and $1$
otherwise.

We verify $j(o(\eta))=P(TM)$ by pairing with surfaces,
following the degree-one argument in
\cite[Proposition~2.24]{DebrayYeYu2026}.
Let $f:\Sigma\to M$ be a smooth map from a closed surface,
possibly nonorientable, transverse to $F$.
Choose disjoint small disks $D_1,\ldots,D_s$ around the
points of $f^{-1}(F)$, and put
$\Sigma_0=\Sigma\setminus\bigcup_i\intt D_i$.
The structure $f^*\eta$ is defined on $(f^*TM)|_{\Sigma_0}$.
Let $\varepsilon_i\in\Ftwo$ be the obstruction to extending
its restriction from $\partial D_i$ over $D_i$. That is, writing $\iota_i:(D_i,\partial D_i)\hookrightarrow(\Sigma,\Sigma_0)$
for the inclusion of pairs, we define
\[\varepsilon_i:=\langle(f\circ\iota_i)^*o(\eta),[D_i,\partial D_i]\rangle\in\Ftwo.\]
By excision and naturality of the relative obstruction,
$\displaystyle\langle j(o(\eta)),f_*[\Sigma]\rangle=\sum_i\varepsilon_i$.

Choose a cell decomposition of $\Sigma$ in which each
$D_i$ is a $2$-cell and $\Sigma_0$ is a subcomplex.
The structure $f^*\eta$ gives a lift over the $1$-skeleton.
The resulting obstruction cocycle is zero on the
$2$-cells in $\Sigma_0$ and takes the value $\varepsilon_i$
on $D_i$.
The universal lifting obstruction is $w_2$ for Spin and
$\Pinp$, and $w_2+w_1^2$ for $\Pinm$
\cite[\S1]{KirbyTaylor1990}.
Thus this cocycle represents $f^*P(TM)$, and
\[
\bigl\langle P(TM),f_*[\Sigma]\bigr\rangle
=\sum_i\varepsilon_i
=\bigl\langle j(o(\eta)),f_*[\Sigma]\bigr\rangle.
\]
Every class in $H_2(M;\Ftwo)$ is represented by a map
from a closed surface, which may be made transverse to $F$.
Hence $j(o(\eta))=P(TM)$.

The absence of higher obstructions also shows that
$\eta$ extends across a component of $F$ exactly when
the corresponding value of $o(\eta)$ is zero.
Therefore, nonextension across every component means
$o(\eta)=U_F$, and implies $P(TM)=j(U_F)$.

Conversely, suppose $P(TM)=j(U_F)$.
Its restriction to $A$ vanishes, so a $\psi$-structure
$\eta_0$ exists on $A$, refining the specified
$\xi$-structure.
By exactness of
\[
H^1(A)\xrightarrow{\delta}H^2(M,A)
\xrightarrow{j}H^2(M),
\]
there is an $x\in H^1(A)$ with
$\delta x=U_F-o(\eta_0)$.
Changing $\eta_0$ by $x$ changes its relative obstruction
by $\delta x$: a cocycle representing $x$ changes the
lifting data on $A$, and the coboundary of an extension
of this cochain represents the change in the obstruction.
Consequently, $o(\eta_0+x)=o(\eta_0)+\delta x=U_F$.
The resulting structure is nonextendable across every
component of $F$.
\end{proof}
The boundary version of \cite[Theorem~2.4]{KirbyTaylor1990}
also gives $[G]=\PD_W(P(TW))$ in
$H_{n-1}(W,\partial W;\mathbb Z/2)$ for a bordism in
\cref{def:char}; the specified complement structures at its
ends remain part of the bordism data.

\begin{theorem}\label{thm:charfib}
For a triple in \eqref{eq:triples}, deletion and filling define
mutually inverse isomorphisms
\[
\begin{aligned}
 d:\Omega_n^{(\xi,P)}
 &\longrightarrow\Omega_n^{\mathrm{nb},\psi\text{-}\mathrm{fib}},\\
 [M,F,\eta]&\longmapsto
 [M\setminus\intt D(\nu_{F/M}),F,
   \pi_{\nu_{F/M}},\nu_{F/M},\id],\\[3pt]
 f:\Omega_n^{\mathrm{nb},\psi\text{-}\mathrm{fib}}
 &\longrightarrow\Omega_n^{(\xi,P)},\\
 [X,B,p,V,\varphi]&\longmapsto
 [X\cup_\varphi D(V),B,\eta_M].
\end{aligned}
\]
For $d$, the complement carries $\eta|_X$ and its boundary is
identified with $S(\nu_{F/M})$. For $f$, the defect is the zero
section, and $\eta_M$ extends the structure of $X$ fiberwise over
$D(V)\setminus B$, matching the boundary collar.
\end{theorem}
\begin{proof}
\emph{The deletion map $d$.}
The nonextension condition on $\eta$ makes every normal
$S^1$-fiber non-$D^2$-bounding, so the displayed tuple is a
cycle of \cref{def:nb}. Tubular neighborhoods are unique up to
isotopy; deleting a tubular neighborhood along such an isotopy
gives a fibered bordism between the resulting cycles, including
their projections and bundle identifications.
For a characteristic bordism $(W,G,\eta_W)$, choose a tubular
neighborhood of $G$ that is product near $\partial W$.
Its complement has ordinary faces equal to the deleted end
cycles and fibered face $S(\nu_{G/W})\to G$, with the identity
sphere-bundle identification. Nonextension across each component
of $G$ gives the required fiber condition throughout this face.
Thus $d$ is well-defined on bordism classes.

\emph{The filling map $f$.}
Boundary collars make $M=X\cup_\varphi D(V)$ a smooth manifold.
The identity $p\circ\varphi=\pi_V$ ensures that $D(V_b)$
fills the specified fiber $p^{-1}(b)$.
In the $O$ cases no ambient orientation is required.
In the $SO$ case the Spin structure orients $X$ and its boundary.
Writing $\pi:S(V)\to B$, the tangent bundle splitting gives
\[
 0=w_1(TS(V))=\pi^*w_1(TB\oplus V)
 =\pi^*(w_1(TB)+w_1(V)).
\]
The mod-two Gysin sequence makes
$\pi^*:H^1(B;\mathbb Z/2)\to H^1(S(V);\mathbb Z/2)$
injective. Thus $TB\oplus V$, and hence $D(V)$, is orientable.
Since the fibers are connected, there is a unique orientation
on each component of $D(V)$ whose boundary orientation is
opposite to that transported from $\partial X$.
These orientations glue to an orientation on $M$.

The extension $\eta_M$ in the statement refines this
ambient structure. Since $S(V)\hookrightarrow D(V)\setminus B$
is a homotopy equivalence, this extension is unique up to
homotopy relative to its boundary value.
Its normal fibers retain the prescribed non-$D^2$-bounding
structure, so it is nonextendable across every component of $B$.
Hence the filling is a characteristic cycle; by
\cref{thm:extension}, $[B]=\PD_M(P(TM))$.

For bordism invariance, take a bordism $(W,N,q,\mathcal V,\Phi)$
in \cref{def:nb} with end cycles $(X_i,B_i,p_i,V_i,\varphi_i)$,
$i=0,1$. Attach $D(\mathcal V)$ along $\Phi$ and
smooth the corners to obtain a manifold $\widehat W$ with
neat defect $N$. The restrictions of the boundary data are $\varphi_i$ and $p_i$.
Thus its ends are the filled cycles $(M_i,B_i,\eta_{M_i})$. In the $SO$ case the preceding orientation
argument applies over $N$; uniqueness makes the orientations
at the ends agree with those used to define $f$.
Fiberwise extension gives a complement structure on
$\widehat W\setminus N$ restricting to the specified $\eta_{M_i}$.
The fiber condition on $S(\mathcal V)$ makes it nonextendable
across every component of $N$. Thus the filled bordism satisfies
\cref{def:char}, and $f$ is well-defined.

\emph{The two compositions.}
Filling a deleted tubular neighborhood restores $M$, $F$,
and the ambient structure. It also restores $\eta$ up to
homotopy, by the uniqueness of extension over the punctured
tube. Hence $f\circ d=\id$.
Conversely, deleting the attached disk bundle restores $X$
with its structure. The boundary identification $\varphi$
recovers the original bundle data, including
$p=\pi_V\circ\varphi^{-1}$, so $d\circ f=\id$.
Both maps preserve disjoint union and are therefore inverse
group isomorphisms.
\end{proof}

\subsection{Structures on an \texorpdfstring{$S^1$}{S1}-bundle}

Let \(V \to B\) be a smooth real $2$-plane bundle equipped with a metric,
and set \(\sigma:=\det V\).
Let \(\pi \colon S(V) \to B\) be its unit sphere bundle.
Write \(T_\pi S(V):=\ker(d\pi \colon TS(V)\to\pi^*TB)\).
Choosing a metric on \(TS(V)\), we obtain
\(T_\pi S(V)^\perp \cong TS(V)/T_\pi S(V)\cong\pi^*TB\), and hence
\[
TS(V)\cong\pi^*TB\oplus T_\pi S(V).
\]

We show that \(T_\pi S(V)\cong\pi^*\sigma\).
Fix \(b\in B\) and \(u\in S(V_b)\), and use the inclusion
\(\iota\colon S(V_b)\hookrightarrow S(V)\) to regard
\(T_uS(V_b)\) as a subspace of \(T_uS(V)\).
Since \(\pi\circ\iota\) is constant, \(d\pi_u\circ d\iota_u=0\),
so \(T_uS(V_b)\subseteq\ker(d\pi_u)\).
Moreover, \(\pi\) is a submersion, and therefore
\(\dim\ker(d\pi_u)=\dim T_uS(V)-\dim T_bB=\dim T_uS(V_b)\).
Thus \((T_\pi S(V))_u=T_uS(V_b)\).
Since \(S(V_b)\) is the unit sphere in \(V_b\), its tangent space at \(u\)
is naturally identified with \(u^\perp\subset V_b\).
Under this identification, define
\[
\varphi\colon T_\pi S(V)\longrightarrow\pi^*\sigma,
\qquad
(u,v)\longmapsto(u,u\wedge v).
\]
For each \(u\), the map \(v\mapsto u\wedge v\) is an isomorphism
from \(u^\perp\) to \(\sigma_b=\det V_b\): it is injective because \(v\perp u\),
and both spaces are one-dimensional.
The map \(\varphi\) is smooth, as seen in local orthonormal trivializations.
It is therefore a smooth bundle isomorphism, giving
\begin{equation}\label{eq:tangent}
TS(V)\cong\pi^*(TB\oplus\sigma).
\end{equation}

\begin{lemma}[$S^1$-bundle structure correspondence]\label{lem:circle}
For $\psi\in\{\Spin,\Pinm,\Pinp\}$, pullback using
\eqref{eq:tangent} induces a bijection of isomorphism classes
\[
 \psi(TB\oplus\sigma)\xrightarrow{\cong}
 \psi(TS(V))_{\mathrm{nb}}.
\]
It is compatible with bundle isomorphisms covering diffeomorphisms
and with restriction to boundaries using product collars.
\end{lemma}
\begin{proof}
All cohomology groups have coefficients in $\Ftwo$, and we write
$S^1=\mathbb R/\mathbb Z$.
For a metric bundle $E$ of rank $k$, let $\operatorname{Fr}(E)$
be its orthonormal frame bundle, using oriented frames for Spin.
A $\psi$-structure $\zeta$ determines an equivariant double cover
$q_\zeta:P_\zeta\to\operatorname{Fr}(E)$ for
$\rho:G\to H$, where $(G,H)=(\Spin_k,SO_k)$ or
$(\Pinpm_k,O_k)$.
A loop $\mathcal F:S^1\to\operatorname{Fr}(E)$ has a closed lift
if it lifts to a continuous map $S^1\to P_\zeta$.

\emph{The local type of a pullback.}
Let $\theta$ be a $\psi$-structure on $TB\oplus\sigma$ and
put $\eta=\pi^*\theta$ using \eqref{eq:tangent}.
Fix $b\in B$, an orthonormal basis $e$ of $T_bB$, and an
orthonormal basis $(v_1,v_2)$ of $V_b$.
For Spin, require $(e,v_1\wedge v_2)$ to have the orientation
specified by $\theta$.
For $\alpha\in S^1$, set
\[
u(\alpha)=\cos 2\pi\alpha\,v_1+\sin 2\pi\alpha\,v_2,
\qquad
t(\alpha)=-\sin 2\pi\alpha\,v_1+\cos 2\pi\alpha\,v_2.
\]
In the chosen splitting along $S(V_b)$, the frames
$(e,t(\alpha))$ define a loop
\[
\mathcal F:S^1\longrightarrow
\operatorname{Fr}\bigl(TS(V)|_{S(V_b)}\bigr).
\]
Since $u(\alpha)\wedge t(\alpha)=v_1\wedge v_2$,
\eqref{eq:tangent} sends these frames to the fixed frame
$(e,v_1\wedge v_2)$ over $b$.
Pulling back a lift of this fixed frame in $P_\theta$
along $\pi|_{S(V_b)}:S(V_b)\to\{b\}$ gives a closed lift
of $\mathcal F$ in $P_\eta$.

Let $\varepsilon$ be the unit vector in the normal line added
for boundary stabilization. With the convention above,
the outward-normal identification is given, for $\lambda\in\mathbb R$, by
\[
\Phi_b:\underline{T_bB}\oplus\underline{\mathbb R}\oplus TS(V_b)
\xrightarrow{\cong}(\underline{T_bB}\oplus TD(V_b))|_{S(V_b)},
\quad (w,\lambda\varepsilon,v)\longmapsto(w,\lambda u+v).
\]
Stabilization sends the closed lift of $(e,t(\alpha))$ to a closed lift of $(e,\varepsilon,t(\alpha))$, and $\Phi_b$ sends the latter frame to $(e,u(\alpha),t(\alpha))$.

Suppose, for a contradiction, that $\eta$ is $D^2$-bounding
on this fiber, and choose a disk structure $\zeta$ on
$\underline{T_bB}\oplus TD(V_b)$ inducing it under $\Phi_b$.
The constant frame $(e,v_1,v_2)$ lifts over the contractible
disk to $\widetilde s:D(V_b)\to P_\zeta$.
With $j_b:S(V_b)\hookrightarrow D(V_b)$ the inclusion, put
$p=\widetilde s\circ j_b\circ u:S^1\to P_\zeta$;
then $p(1)=p(0)$.
Put $m=\dim B$. Relative to the constant frame, the rotating
frame is given by
\[
R:S^1\longrightarrow SO_{m+2},\qquad
R(\alpha)=\operatorname{diag}\!\left(I_m,
\begin{pmatrix}
\cos 2\pi\alpha&-\sin 2\pi\alpha\\
\sin 2\pi\alpha&\cos 2\pi\alpha
\end{pmatrix}\right).
\]
This loop represents a generator of $\pi_1(SO(m+2))$.
By \cite[\S1, p.~180]{KirbyTaylor1990} and homotopy lifting,
its lift $\widetilde R:[0,1]\to G$ starting at $1$ ends at
$-1$, the nonidentity element of $\ker\rho$.
Here and below, loops are evaluated on $[0,1]$ via its
quotient map to $S^1$.
Equivariance gives
\[
q_\zeta\bigl(p(\alpha)\cdot\widetilde R(\alpha)\bigr)
=q_\zeta(p(\alpha))\cdot\rho(\widetilde R(\alpha))
=(e,u(\alpha),t(\alpha)).
\]
This lift ends at $p(0)\cdot(-1)\ne p(0)$.
By uniqueness of path lifting, the other lift is its translate
by $-1$, so neither is closed.
This contradicts the closed lift in the stabilized $\eta$,
since an isomorphism of structures over the identity preserves
closed lifts. Hence $\eta$ is fiberwise non-$D^2$-bounding.

\emph{Existence on the base.}
Conversely, let $\eta$ be fiberwise non-$D^2$-bounding.
Let $P$ be the obstruction class in \eqref{eq:triples}:
$P=w_2$ for Spin and $\Pinp$, and $P=w_2+w_1^2$ for $\Pinm$.
The mod-two Gysin sequence gives
\begin{equation}\label{eq:gysin}
0\longrightarrow H^1(B)
\xrightarrow{\pi^*}H^1(S(V))
\xrightarrow{\pi_!}H^0(B)
\xrightarrow{\smile w_2(V)}H^2(B),
\end{equation}
where $\pi_!$ evaluates on each fiber.
For Spin, \eqref{eq:tangent} and injectivity of $\pi^*$ imply
$w_1(TB)+w_1(V)=0$.
Connected fibers then give a unique orientation on
$TB\oplus\sigma$ inducing that of $\eta$.
The orientation argument in \cref{thm:charfib} also gives an
orientation on $D(V)$ inducing that of $\eta$ on its boundary.

Using $D(V)\setminus B\cong S(V)\times(0,1]$,
fiberwise extension of $\eta$ by the product structure defines
$\zeta$ on $TD(V)|_{D(V)\setminus B}$.
As in the relative obstruction argument of
\cref{thm:extension}, its extension obstruction lies in $o(\zeta)\in H^2(D(V),D(V)\setminus B)\cong H^0(B)$.
Its restriction to each normal disk is $1$, precisely because
the specified fiber structure is non-$D^2$-bounding.
Thus $o(\zeta)$ is the Thom class $U_V$.
Its absolute image is $P(TD(V))$; pulling back to the zero
section gives
\[
P(TB\oplus V)=w_2(V),
\]
using the zero-section pullback of the Thom class
(cf.\ \cite[Corollary~2.23]{DebrayYeYu2026}).
The Whitney formulas and $\sigma=\det V$ give
\[
w_1(TB\oplus V)=w_1(TB\oplus\sigma),\qquad
P(TB\oplus V)=P(TB\oplus\sigma)+w_2(V).
\]
Consequently $P(TB\oplus\sigma)=0$, so a base structure
exists, with the orientation already specified for Spin.

\emph{Uniqueness and surjectivity.}
Choose a base structure $\theta_0$ as above.
By the $H^1$-torsor description
\cite[\S1]{KirbyTaylor1990}, write
$\eta=\pi^*\theta_0+x$ with $x\in H^1(S(V))$.
Both structures restrict to the unique nonbounding lifting
type on each fiber, so $\pi_!x=0$.
Exactness of \eqref{eq:gysin} gives a unique $y\in H^1(B)$
with $x=\pi^*y$, hence $\eta=\pi^*(\theta_0+y)$.
This proves the bijection for each Spin orientation;
distinct orientations remain distinct under pullback.

\emph{Compatibility.}
For a base with boundary, take the splitting
$TS(V)\cong\pi^*TB\oplus T_\pi S(V)$ to be a product on
boundary collars. Then \eqref{eq:tangent} sends the outward
normal of $S(V)$ to the pullback of the outward normal of $B$;
removing these normal lines gives the boundary correspondence.
For a bundle isomorphism $A:V\to V'$ covering a diffeomorphism,
transport the metrics and this splitting.
The identity $\det A(u\wedge v)=Au\wedge Av$ shows that
\eqref{eq:tangent} commutes with the induced bundle maps.
Metrics and splittings are unique up to homotopy, so this
proves compatibility on isomorphism classes of structures.
The inverse correspondence is compatible by uniqueness.
\end{proof}

\begin{remark}[Thom-spectrum convention]\label{rem:thom}
We use the twisted bordism notation of \cite{DebrayYeYu2026}:
\[
 \Omega_k^\psi((BO_2)^{\sigma-1})
 =\pi_k\bigl(MT\psi\wedge(BO_2)^{\sigma-1}\bigr).
\]
For a genuine bundle $E\to X$, $\Th(X,E)$ denotes its based
Thom space.  Here $MT\psi$ represents stable tangential $\psi$-bordism,
$BO_2$ is the unbased classifying space, and $(BO_2)^{\sigma-1}$
is the Thom spectrum of the virtual rank-zero bundle $\sigma-1$.
Geometrically, cycles are $(B^k,V,\theta)$ with $V\to B$ of rank
two and $\theta$ a stable $\psi$-structure on $TB\oplus\det V$;
bordisms carry the same data extending those at their boundaries.
\end{remark}

\begin{corollary}\label{cor:twist}
For $\psi=\Spin,\Pinm,\Pinp$ there is an isomorphism
\[
 \Omega_k^{\mathrm{nb},\psi\text{-}\mathrm{sph}}
 \cong\Omega_k^\psi((BO_2)^{\sigma-1}).
\]
\end{corollary}
\begin{proof}
Send $(B,V,\eta)$ to $(B,V,\theta)$, where $\theta$ is the unique
preimage under \cref{lem:circle}.  The same lemma applies on a
compact base bordism and commutes with its boundary structures.
It therefore sends a sphere-bundle bordism satisfying
\cref{def:nb} to a twisted bordism.  Conversely, pullback sends
a twisted bordism to a sphere-bundle bordism satisfying that
fiber condition.  These operations are inverse on
cycles and on bordism classes; they preserve disjoint union.
\end{proof}

\subsection{The characteristic long exact sequence}

\begin{theorem}[Main theorem]\label{thm:main}
Fix one of the three triples in \eqref{eq:triples}, and define
$\Omega_n^{(\xi,P)}$ as in \cref{def:char}.
Let $\sigma\to BO_2$ be the universal determinant line.
There is a long exact sequence
\[
 \cdots\longrightarrow
 \Omega_{n-1}^{\psi}((BO_2)^{\sigma-1})
 \xrightarrow{t}\Omega_n^\psi
 \xrightarrow{i}\Omega_n^{(\xi,P)}
 \xrightarrow{R}\Omega_{n-2}^{\psi}((BO_2)^{\sigma-1})
 \xrightarrow{t}\Omega_{n-1}^\psi\longrightarrow\cdots.
\]
\end{theorem}
\begin{proof}
Apply \cref{prop:nb}, replace its fibered-boundary group using
\cref{thm:charfib}, and replace its sphere-bundle groups using
\cref{cor:twist}. Transporting the maps along these isomorphisms
gives the required exact sequence.
\end{proof}

The map $i$ sends $[M,\eta]$ to $[M,\varnothing,\eta]$, forgetting
to $\xi$ on the ambient manifold.  The map
\[
 R([M,F,\eta])=[F,\nu_{F/M},\theta_\eta]
\]
uses the boundary structure of the deleted complement and the
inverse correspondence of \cref{lem:circle}.  Finally,
$t([B,V,\theta])=[S(V),\pi^*\theta]$, using \eqref{eq:tangent}.
The induced maps on bordism groups do not depend on these choices.

\paragraph{Names of the three theories.}
In the following corollaries the names GM and $\mathrm{KT}^\pm$
refer to the groups of \cref{def:char}.

\begin{corollary}[Guillou--Marin type]\label{cor:gm}
For $(\xi,P,\psi)=(SO,w_2,\Spin)$, write
$\Omega_n^{\mathrm{GM}}:=\Omega_n^{(SO,w_2)}$ with the preceding
convention.  Then
\[
 \cdots\to\Omega_{n-1}^{\Spin}((BO_2)^{\sigma-1})
 \to\Omega_n^{\Spin}\to\Omega_n^{\mathrm{GM}}
 \to\Omega_{n-2}^{\Spin}((BO_2)^{\sigma-1})
 \to\Omega_{n-1}^{\Spin}\to\cdots.
\]
\end{corollary}

\begin{corollary}[Kirby--Taylor$^-$ type]\label{cor:ktm}
For $(\xi,P,\psi)=(O,w_2+w_1^2,\Pinm)$, write
$\Omega_n^{\mathrm{KT}^-}:=\Omega_n^{(O,w_2+w_1^2)}$ with the
preceding convention.  Then
\[
 \cdots\to\Omega_{n-1}^{\Pinm}((BO_2)^{\sigma-1})
 \to\Omega_n^{\Pinm}\to\Omega_n^{\mathrm{KT}^-}
 \to\Omega_{n-2}^{\Pinm}((BO_2)^{\sigma-1})
 \to\Omega_{n-1}^{\Pinm}\to\cdots.
\]
\end{corollary}

\begin{corollary}[Kirby--Taylor$^+$ type]\label{cor:ktp}
For $(\xi,P,\psi)=(O,w_2,\Pinp)$, write
$\Omega_n^{\mathrm{KT}^+}:=\Omega_n^{(O,w_2)}$ with the preceding
convention.  Then
\[
 \cdots\to\Omega_{n-1}^{\Pinp}((BO_2)^{\sigma-1})
 \to\Omega_n^{\Pinp}\to\Omega_n^{\mathrm{KT}^+}
 \to\Omega_{n-2}^{\Pinp}((BO_2)^{\sigma-1})
 \to\Omega_{n-1}^{\Pinp}\to\cdots.
\]
\end{corollary}

\end{document}